\documentclass[11pt]{article}
\usepackage{graphicx}
\usepackage{amsmath,amssymb,amsthm,amsfonts
}
\usepackage{amssymb}
\newtheorem{thm}{Theorem}[section]
\newtheorem{cor}{Corollary}[section]
\newtheorem{lem}{Lemma}[section]

\newtheorem{rem}{Remark}[section]
\theoremstyle{definition}

\numberwithin{equation}{section}

\DeclareMathSymbol{\C}{\mathalpha}{AMSb}{"43}

\newcommand{\eps}{\varepsilon}

\newcommand{\R}{{\mathbb{R}}}
\newcommand{\h}{{\mathcal{H}}}
\newcommand{\inte}{\int_{\mathbb{R}^2}}

\newcommand{\bsub}{\begin{subequations}}
\newcommand{\esub}{\end{subequations}$\!$}

\begin{document}

\title{Properties of ground states of   attractive Gross-Pitaevskii equations with multi-well  potentials\thanks{Email: YJG:yjguo@wipm.ac.cn; ZQW:zhi-qiang.wang@usu.edu;  XYZ:xyzeng@wipm.ac.cn; HSZ:hszhou@wipm.ac.cn; Tel:+86-2787199196}}

\author{Yujin Guo$^a,^b$,\
    Zhi-Qiang Wang$^c$,\ Xiaoyu Zeng$^a,^b$\ \ and Huan-Song Zhou$^a,^b$\\
\small \it $^a$Wuhan Institute of Physics and Mathematics,
    Chinese Academy of Sciences,\\
    \small \it P.O. Box 71010, Wuhan 430071, P. R. China\\
    \small \it $^b$ National Center for Mathematics and Interdisciplinary Sciences, Beijing 100190, P.R. China\\
   \small \it $^c$Department of Mathematics, Utah State University, Logan, UT 84322, USA \\
  }

\date{\today}

\smallbreak \maketitle

\begin{abstract}
We are interested in the attractive  Gross-Pitaevskii (GP) equation in $\R^2$, where  the external potential $V(x)$ vanishes on $m$ disjoint bounded domains $\Omega_i\subset \R^2\  (i=1,2,\cdots,m)$ and $V(x)\to\infty$ as $|x|\to\infty$, that is, the union of these $\Omega_i$ is the bottom of the potential well. By making some delicate estimates on the energy functional of the GP equation, we prove that when the interaction strength $a$ approaches some critical value $a^*$ the ground states concentrate and blow up at the center of the incircle of some $\Omega_j$ which has the largest inradius. Moreover, under some further conditions on $V(x)$ we show that the ground states of GP equations are unique and radially symmetric at leat for almost every $a \in (0, a^*)$.

\end{abstract}

\vskip 0.2truein


\vskip 0.2truein

\section{Introduction}

In this paper, we consider the following stationary (i.e. time-independent) Gross-Pitaevskii (GP) equation in $\R^2$:
\begin{equation}\label{d1.1}
-\Delta u+V(x)u=\mu u+a |u|^2u, \quad  x\in
\mathbb{R}^2,
\end{equation}
where $V(x)$ is the external potential, $\mu\in\R$ is the chemical potential, $a>0$ represents the attractive  interaction strength. (\ref{d1.1}) is a model equation for the single-particle wave function in a Bose-Einstein condensate (BEC, in short).

 It is well-known that equation (\ref{d1.1}) can be obtained from the time dependent GP equation when we look for the standing wave type solutions $\psi(x,t)=u(x)e^{-i\mu t}$, where $i^2=-1$. Moreover, (\ref{d1.1}) is also the Euler-Lagrange equation of the following constrained minimization problem
 \begin{equation}\label{eq1.3}
e(a):=\inf_{u\in \mathcal{M}} E_a(u)\,,
\end{equation}
where $E_a(u)$ is the so-called GP functional \cite{D,P}
\begin{equation}\label{f}
  E_a(u):=\int_{\R ^2}\Big(|\nabla
  u|^2+V(x)|u|^2\Big)dx-\frac{a}{2}\inte|u|^4dx, \  u\in\h.
\end{equation}
Here we define
\begin{equation}
   \h := \Big \{u\in  H^1(\R ^2):\ \int _{\R ^2}  V(x)|u(x)|^2 dx<\infty \Big\}, \label{H}
\end{equation}
and  \begin{equation}\label{norm}\mathcal{M}=\Big\{ u\in\mathcal{H}; \inte|u|^2dx=1\Big\}\,.
\end{equation}

Very recently, by assuming that
 \begin{itemize}
\item[$(V)$]  $0\leq V(x)\in L_{loc}^\infty(\R^2)$,
$\underset{|x|\to\infty}{\lim} V(x) = \infty$,
\end{itemize}
the existence and non-existence of minimizers for (\ref{eq1.3}) are proved in \cite{Bao2,GS}, which show  that there exists a critical value $a^*>0$ such that (\ref{eq1.3}) has at least one minimizer if $a\in [0,a^*)$, and (\ref{eq1.3}) has no minimizers if $a\geq a^*$. Furthermore, the critical value $a^*=\|Q(x)\|_2^2$, that is, the square of the $L^2$-norm of the unique positive solution of the famous nonlinear scalar field equation
\begin{equation}
-\Delta u+ u-u^3=0\  \mbox{  in } \  \R^2,\ \  u\in
H^1(\R ^2),  \label{Kwong}
\end{equation}
see, e.g.  \cite{GNN,K,Li,mcleod}.

 We mention that the above conclusions give a rigid mathematical explanation for the collapse of attractive BEC ($a>0$), that is, if the particle number increases  beyond a critical value $a^*$ in an attractive BEC system, the collapse must occur, see \cite{D,HM,KM}, etc. In this paper, we aim to investigate the details of this collapse under certain general potentials. Roughly speaking, we analyze how a ground state, that is, a least energy solution of (\ref{d1.1}), blows up as $a\nearrow a^*$. So, we introduce now a rigorous definition of the ground states of (\ref{d1.1}). With $\mu$ fixed we first define the associated energy functional of equation (\ref{d1.1}) by
 \begin{equation}\label{3.1.5}
J_{a,\mu}(u):=\frac{1}{2}\inte \Big(|\nabla u|^2+\big(V(x)-\mu\big)|u|^2\Big)dx-\frac{a}{4}\inte |u|^4dx,\ \ u\in\h .
\end{equation}
 Then, $u\in\h\setminus \{0\}$ is called a nontrivial (weak) solution of (\ref{d1.1}) if $\langle J_{a,\mu}'(u), \varphi\rangle=0$ for any $\varphi\in\h$. Let
 \begin{equation}\label{eq1.8}S_{a,\mu}:=\big\{u(x): \text{  $u(x)$ is a nontrivial  solution of (\ref{d1.1})} \big\},\end{equation}
 and define
 \begin{equation}\label{3.1.7}
G_{a,\mu}:=\big\{u(x)\in S_{a,\mu}:\   J_{a,\mu}(u)\leq J_{a,\mu}(v) \text{ for all $v\in S_{a,\mu}$} \big\}.
\end{equation}
Therefore, we say $u\in\h$ is a {\it ground state} (or, a {\it least energy solution}) of (\ref{d1.1}) if $u\in G_{a,\mu}$. As the least energy solutions are one sign solutions we will always assume a ground state is a positive solution. Now, a natural question is whether the minimizers of (\ref{eq1.3}) are ground states of (\ref{d1.1}) and how about the converse?  Our first theorem is to answer these questions. For any $a\in [0, a^*)$, let
\begin{equation}\label{eq1.10}
\Lambda_a:=\big\{u_a: u_a \text{ is a minimizer of $e(a)$ in (\ref{eq1.3})}\big\}.
\end{equation}
If $u_a\in \Lambda_a$, as illustrated in \cite{GS}, we may assume that $u_a\geq0$ and $u_a$ satisfies (\ref{d1.1}) with a suitable $\mu=\mu_a$.
By the results of \cite{GZZ} and the references therein, we known that $e(a)$ defined in (\ref{eq1.3}) has a unique positive minimizer for any $a>0$ being small enough ($a<a^*$). So, we may define
\begin{equation}\label{1.15}
a_{*}:=\sup\big\{l>0: \text{ $e(a)$ has a unique positive minimizer for all $a\in[0,l)$} \big\},
\end{equation}
 and $0<a_*\leq a^*$. Our Theorem \ref{pro1} shows that $\Lambda_a=G_{a,\mu_a}$ for some $\mu_a\in\R$.
\begin{thm}\label{pro1}
Let the condition $(V)$ be satisfied. Then, for all $a\in [0,a_*)$ and for a.e. $a\in[a_{*},a^*)$, all minimizers of $e(a)$ satisfy  (\ref{d1.1}) with a fixed    Lagrange multiplier $\mu=\mu_a$  and $\Lambda_a =G_{a,\mu_a}$.
\end{thm}

We remark that similar conclusions to Theorem \ref{pro1} was also studied elsewhere for different types of problems, see for instance \cite[Chapter 8]{C}. Theorem \ref{pro1} is proved in Section 2, where some  fundamental properties of minimizers and the energy $e(a)$ are also  addressed. It expects that the fact $\Lambda_a =G_{a, \mu_a}$ is useful for the  further understanding of the minimizers for (\ref{eq1.3}) as well as  ground states of (\ref{d1.1}). Based on Theorem \ref{pro1}, we have  the following result on the  uniqueness of positive minimizers of (\ref{eq1.3}) under some further condition of  potential $V(x)$.

\begin{cor}\label{pro1.4}
Suppose $V(x)$ satisfies  condition $(V)$ as well as
\begin{equation}\label{1.17}
V(x)=V(|x|)\ \text{ and }\ V'(|x|)\geq 0(\not\equiv0).
\end{equation}
Then, for all $a\in[0,a_{*})$ and a.e. $a\in[a_{*},a^*)$, $e(a)$ has a unique positive minimizer which must be radially symmetric about the origin.
\end{cor}
In view of Theorem \ref{pro1}, we start to study the properties of the ground states for (\ref{d1.1}) by investigating those  of the minimizers for  (\ref{eq1.3}). If $V(x)$ satisfies $(V)$ and has finitely many isolated zero points, e.g. $V(x)=\overset{n}{\underset{i=1}\prod}
|x-x_i|^{p_i}$ with $p_i>0$ and  $x_i\not =x_j$ for $i\not =j$, the detailed analysis for the concentration and symmetry breaking of the minimizers of (\ref{eq1.3}) as $a\nearrow a^*$ was first studied in \cite{GS} based on some precise estimates of the GP energy $e(a)$. However, the methods  used in \cite{GS} depend heavily on the potential $V(x)$ having a finite number of zeros $\{x_i\in\R^2: i=1,2,\cdots,n\}$. Very recently, the results of \cite{GS} were extended to the case where
\begin{equation}\label{1.16}
 V(x)=(|x|-A)^2,\  A>0 \text{ and } x\in\R^2,
\end{equation}
see \cite{GZZ} for the details. Clearly, the above $V(x)$ has infinitely many zeros, that is, $\{x\in\R^2: |x|=A\}$, which has zero measure. Then, the way of \cite{GS} for getting the optimal energy estimates for the GP functional (\ref{f}) does not work anymore. In \cite{GZZ}, the authors provided a new approach for establishing the energy estimates under the potential (\ref{1.16}). As a continuation of \cite{GS,GZZ}, in the present paper we want to consider problem (\ref{eq1.3}) for more general potential $V(x)$ which may not have an explicit expression  like (\ref{1.16}), etc. Particularly, we allow $V(x)$ vanishes on a set with nonzero measure. These new features on $V(x)$ cause some essential difficulties  on the estimates of GP energy. If the condition $(V)$ is slightly strengthened, we have the following general theorem  on the concentration behavior of the minimizers of (\ref{eq1.3}).
\begin{thm}\label{thm1.2} Let $u_a\geq0$ be
a nonnegative minimizer of (\ref{eq1.3}) for $a\in(0,a^*)$. If $V(x)$ satisfies
\begin{itemize}
\item[$(V_1)$]  $V\in C^\alpha_{loc}(\R^2)$ with some $\alpha\in(0,1)$, and
$\underset{|x|\to\infty}{\lim} V(x) = \infty$ as well as $\underset{x\in \R^2}{\inf}
V(x) =0$,
\end{itemize}
then, for  any sequence
$\{a_k\}$ with $a_k\nearrow a^*$ as $k\to\infty$, there exists a
subsequence, still denoted by $\{a_k\}$, of $\{a_k\}$ such that each $u_{a_k}$ has a
unique maximum point $\bar z_k$, which satisfies $\underset{k\to\infty}\lim\bar z_k= x_0$ and $V(x_0)=0$.
Moreover, we have
\begin{equation}\label{con:a1}
\lim_{k\to\infty}\varepsilon_k u_{a_k}\big(\varepsilon_k x+\bar z_k
\big) =\frac { Q(x)}{\sqrt{a^*}}\ \text{ strongly in }\ H^1(\R^2),
\end{equation}
where $\varepsilon_k>0$ is defined by and satisfies
\begin{equation*}%\label{1:na1}
\varepsilon_k:=\Big(\inte|\nabla u_{a_k}|^2dx\Big)^{-\frac{1}{2}}\to 0\quad \text{ as }\ k\to\infty\,.
\end{equation*}
\end{thm}
We note that in Theorem \ref{thm1.2}, we could not give explicitly  the convergent rates for $\eps_k$ and $u_{a_k}$ as $k\to\infty$.
When the potential $V(x)$ is either polynomial as in \cite{GS} or ring-shaped as in \cite{GZZ}, there are some precise information on the zero points of $V(x)$, from which  we can deduce the exact convergent  rates of $\eps_k$ and $u_{a_k}$. However, the methods of \cite{GS,GZZ} seem to work at most for the case where  the zero set  $\{x\in \R^2:\, V(x)=0\}$ has zero measure. Now,  we address the refined concentration results for the case where the zero set  $\{x\in \R^2:\, V(x)=0\}$  has a positive Lebesgue measure, such that  $V(x)$ is a potential with multiple wells and the bottom of the wells have positive measure.   Towards this purpose, we require some additional conditions on  $V(x)$.
\begin{itemize}
\item[$(V_2)$] $\bar \Omega=\big\{x\in\R^2;\, V(x)=0\big\}$, where $\Omega=\overset{m}{\underset{i=1}\bigcup}\Omega_i$  and  $\Omega _1,\cdots , \Omega _m$ are disjoint bounded domains in $\R^2$, and $V(x)>0$ in $\bar \Omega ^c:=\R^2\setminus \bar \Omega$.
\end{itemize}
Moreover, we denote
\begin{equation}\label{1:R}
R:=\max_{i=1,\cdots, m}R_i,\ \text{ where }\ R_i:=\max_{x\in\Omega_i}dist (x,\partial \Omega_i)>0\,,
\end{equation}
and \begin{equation}\label{1.18}
\Lambda:=\big\{\Omega_i:\, R_i=R\ \mbox{for some} \ i\in \{1,\cdots,m\}\big\}\,.
\end{equation}
With the further information $(V_2)$ on $V(x)$, we then have the following optimal energy estimates and refined concentration behavior  of nonnegative minimizers of (\ref{eq1.3}) as $a\nearrow a^*$.

\begin{thm}\label{thm1.3}
Suppose  $V(x)$ satisfies  $(V_1)$ and $(V_2)$, then
\begin{itemize}
  \item [\rm(i)] The GP energy $e(a)$ satisfies
  \begin{equation}\label{5.1.4}
e(a)\approx \frac{1}{4R^2a^*}\big(|\ln(a^*-a)|\big)^2(a^*-a)\quad \text{as}\quad a\nearrow a^*\,,
\end{equation}
where $f(a)\approx g(a)$ means that $f/g\to 1$ as $ a\nearrow a^*$.
  \item [\rm(ii)]
For any sequence $\{a_k\}$ satisfying $a_k\nearrow a^*$ as $k\to\infty$, there exists a
subsequence, still denoted by $\{a_k\}$, of $\{a_k\}$ such that each $u_{a_k}$ has a
unique maximum point  $\bar z_k$, which satisfies   $\underset{k\to\infty}{\lim}\bar z_k= x_0\in \Omega_{i_0}$ for some $\Omega_{i_0}\in\Lambda$ and
\begin{equation}\label{1:R0}
\lim_{k\to\infty}dist(\bar z_k,\partial\Omega_{i_0})=R_{i_0}=R\,.
\end{equation}
Moreover, we have
\begin{equation}\label{con:a}
\lim_{k\to\infty}\frac{2R}{|\ln(a^*-a_k)|} u_{a_k}\Big(\bar z_k +
\frac{2R}{|\ln(a^*-a_k)|}x\Big) =
\frac { Q(x)}{\sqrt{a^*}}\text{ strongly in $H^1(\R^2)$}.
\end{equation}

\end{itemize}
\end{thm}

 Theorem \ref{thm1.3} implies that the minimizers of $e(a)$ blow up at the most centered point of $\Omega$ as $a\nearrow a^*$.  In this case, since   the infimum of $V(x)$ attains in the whole domain of $\Omega$, the existing methods of \cite{GS,GZZ} cannot be applied directly, and some different arguments are necessary for proving  Theorem \ref{thm1.3}. Indeed, based on Theorem  \ref{thm1.2} the key point of proving Theorem \ref{thm1.3} is to establish the optimal blow-up rate of $u_{a_k}$ and   a refined description of the unique maximum point  $\bar z_k$.  To achieve these aims, a new and suitable trial function is needed to establish the optimal upper bound for GP energy   $e(a)$, see Lemma \ref{lem2.2} for details. On the other hand, as stated in Lemma \ref{lem2.3}, a proper lower bound  of the minimizers is also necessary in order to analyze the optimal energy  bound of $e(a)$. By a delicate analysis, these results finally yield the optimal blow-up rate of $u_{a_k}$ and a refined description of the unique maximum point  $\bar z_k$. We also remark that the proof of Theorem  \ref{thm1.3}  implies the following refined estimate
\begin{equation}\label{1:na}
\lim_{k\to\infty}\inte|\nabla u_{a_k}|^2dx\Big/\Big(\frac{|\ln(a^*-a_k)|}{2R}\Big)^2=1,
\end{equation}
see  (\ref{3:fin1}) for more details.

Furthermore, Theorem \ref{thm1.3} also indicates  that  {\em symmetry breaking} occurs in the minimizers when the potential $V$ is radially symmetric. For example, suppose that  $V(x)=V(|x|)$ satisfies $(V_1)$ and $V(x)>0$ in $\bar \Omega ^c:=\R^2\setminus \bar \Omega$, where
\begin{equation}\label{1.28}
 \Omega=\big\{ 0<R_1<|x|<R_2: \, V(|x|)=0\big\}
\end{equation}
for positive constants $R_1$ and $R_2$.
It then follows from Theorem \ref{thm1.3} that  all  nonnegative minimizers of (\ref{eq1.3}) can concentrate  at any point on the circle $ |x|=\frac{R_1+R_2}{2} $. This further implies that    there exists an $\bar a$ satisfying
$0< \bar a <a^*$ such that for any $a\in[\bar a , a^*)$, the GP energy $e(a)$
 has infinitely many  different nonnegative minimizers. However, $e(a)$ has a unique nonnegative minimizer $u_a$ for all $a\in[0,a_{*})$, where $a_{*}>0$ is given by (\ref{1.15}), and by rotation $u_a$ must be radially symmetric around the origin. In view of the  above results, we have immediately the following corollary.

\begin{cor}
Suppose $V(x)$ satisfies $(V_1)$ and  (\ref{1.28}). Then   there exist two positive constants $\bar a$ and $a_*$ satisfying $0<a_{*}\le \bar a<a^*$  such that
\begin{itemize}
  \item[(i)]  If $a\in[0,a_{*})$, $e(a)$ has a unique nonnegative minimizer which is radially symmetric about the origin.
  \item [(ii)]  If $a\in[\bar a,a^*)$, $e(a)$ has infinitely many different   nonnegative minimizers which are     not  radially symmetric about the origin.
\end{itemize}
\end{cor}

Finally, we mention  that the  symmetry breaking bifurcation of ground states for nonlinear Schr\"odinger or Gross-Pitaevskii equations has been studied extensively  in the literature, see  e.g.  \cite{Schnee,J,K11,K08}.
Also, the concentration phenomena have also been studied
elsewhere in different contexts, such as \cite{BWang,LW,NW,Wang,WangZ} and the references therein. However, our analysis is mainly involved with either the variational methods or the attractive case, which is different from those employed in  the above mentioned papers.

This paper is organized as follows.  In Section 2 we shall prove Theorem \ref{pro1.4} and Corollary \ref{pro1.4} as well as some other analytical properties of minimizers for  (\ref{eq1.3}). Section \ref{sec3} is devoted to the proof of Theorem \ref{thm1.2} on the  concentration behavior of the minimizers of $e(a)$ under general potentials. In Section \ref{sec4}, for the multiple-well  potentials we  first establish optimal energy estimates of nonnegative minimizers as $a\nearrow a^*$, upon which we then  complete the proof of Theorem \ref{thm1.3} on the refined  concentration behavior of nonnegative minimizers as $a\nearrow a^*$.

\section{Uniqueness of the ground state of (\ref{d1.1})}\label{sec2}
In this section, we first study some properties of the GP energy $e(a)$, upon which we give the proofs for Theorem \ref{pro1.4} and Corollary \ref{pro1.4}. These imply that $\Lambda_a=G_{a, \mu_a}$ and the uniqueness of the ground states of (\ref{d1.1}).

  Before going to the discussion of the properties on $e(a)$, let us recall some auxiliary results which are often used later. By Theorem B of \cite{W}, we have the following Gagliardo-Nirenberg inequality
\begin{equation}\label{GNineq}
\inte |u(x)|^4 dx\le \frac 2 {\|Q\|_2^{2}} \inte |\nabla u(x) |^2dx
\inte |u(x)|^2dx ,\   \  u \in H^1(\R ^2)\,,
\end{equation}
 where the equality is achieved at $u(x) = Q(|x|)$ with $Q(|x|)$ being the unique (up to translations) radial positive solution  of (\ref{Kwong}). Using (\ref{GNineq}) together with (\ref{Kwong}) we know that, see also \cite[Lemma 8.1.2]{C},
\begin{equation}\label{1:id}
\inte |\nabla Q |^2dx  =\inte |Q| ^2dx=\frac{1}{2}\inte |Q| ^4dx.
\end{equation}
 Moreover, it follows  from Proposition 4.1 of  \cite{GNN} that $Q(x)$ satisfies
 \begin{equation}
Q(x) \, , \ |\nabla Q(x)| = O(|x|^{-\frac{1}{2}}e^{-|x|}) \quad
\text{as  $|x|\to \infty$.}  \label{4:exp}
\end{equation}
The following compactness lemma was essentially proved in  \cite[Theorem XIII.67]{RS} or \cite[Theorem 2.1]{BW}, etc.

\begin{lem}\label{2:lem1}
Let  $(V)$ be satisfied. Then, the  embedding
$\h\hookrightarrow L^q(\R^2)$ is compact for all $ q\in[2,\infty)$.
\end{lem}
Now, we give  our result on the smoothness of the GP energy $e(a)$ with respect to $a$.

\begin{lem}  \label{lem2.B}
Let condition $(V)$ be satisfied. Then, for $a\in(0,a^*)$,
the left and right derivatives of $e(a)$ always exist in $[0,a^*)$ and satisfy $$e'_-(a)=-\frac{\alpha_a}{2},\quad e'_+(a)=-\frac{\gamma_a}{2},$$
where
\begin{equation}\label{0.3}
\alpha_a:=\inf\Big\{\inte |u_a|^4dx:\, u_a \in \Lambda_a\Big\},\quad\gamma_a:=\sup\Big\{\inte |u_a|^4dx:\, u_a \in \Lambda_a\Big\},
\end{equation}
and $\Lambda_a$ is given by (\ref{eq1.10}).
\end{lem}

\noindent\textbf{Proof.} Since $V(x)$ satisfies $(V)$, by the definition of $e(a)$, one can derive that $e(a)$ is  decreasing in $a\in[0,a^*)$  and satisfies
$$0\leq \inf_{x\in \R^2} V(x)\leq e(a)\leq e(0)=\mu_1\text{ for all } a\in[0,a^*),$$  where (\ref{GNineq}) is used in the above inequality
 and $\mu_1$ is the first eigenvalue of $-\Delta+V(x)$ in $\h$.
Moreover, it follows from
(\ref{GNineq}) that
\begin{equation}\label{u:bound}\inte |u_a|^4dx\leq \frac{2e(a)}{a^*-a}\leq\frac{2\mu_1}{a^*-a} \quad \text{for}\quad a\in[0,a^*)\,.\end{equation}
For any $a_1,a_2\in[0,a^*)$, we have
\begin{equation}\label{0.6}
e(a_1)\geq e(a_2)+\frac{a_2-a_1}{2}\inte |u_{a_1}|^4dx,\quad \forall\ u_{a_1}\in\Lambda_{a_1},
\end{equation}
\begin{equation}\label{0.7}
e(a_2)\geq e(a_1)+\frac{a_1-a_2}{2}\inte |u_{a_2}|^4dx,\quad \forall\ u_{a_2}\in\Lambda_{a_2},
\end{equation}
and therefore  $\lim_{a_2\to a_1}e(a_2)=e(a_1)$. This implies that
\begin{equation} e(a)\in C([0,a^*),\R^+).\label{u:co} \end{equation}
Furthermore, it follows  from (\ref{0.6}) and (\ref{0.7}) that
\begin{equation}\label{0.8}
\frac{a_2-a_1}{2}\inte |u_{a_1}|^4dx\leq e(a_1)-e(a_2)\leq \frac{a_2-a_1}{2}\inte |u_{a_2}|^4dx.
\end{equation}
Assume $0<a_1<a_2<a^*$. It then follows from  (\ref{0.8}) that
\begin{equation}
 -\frac{1}{2}\inte |u_{a_2}|^4dx\leq \frac{e(a_2)-e(a_1)}{a_2-a_1}\leq  -\frac{1}{2}\inte |u_{a_1}|^4dx, \ \forall\ u_{a_i}\in\Lambda_{a_i},\ i=1,2.
\end{equation}
This implies that
\begin{equation}\label{0.12}-\inf_{u_{a_2}\in\Lambda_{a_2}}\frac{1}{2}\inte |u_{a_2}|^4dx\leq\frac{e(a_2)-e(a_1)}{a_2-a_1}\leq  -\frac{1}{2}\gamma_{a_1}.\end{equation}
In view of (\ref{u:bound}) and Lemma \ref{2:lem1}, there exists $\bar u\in\h$ such that for all $a_2\searrow a_1$,
$$u_{a_2} \rightharpoonup \bar u\text{ weakly in  }\h \text{ and } u_{a_2}\rightarrow \bar u\text{ strongly in  }L^p(\R^2)\text{ for all $p\in[2,\infty)$ }.$$
It then follows from (\ref{u:co}) that
$$e(a_1)=\lim_{a_2\searrow a_1}e(a_2)=\lim_{a_2\searrow a_1}E_{a_2}(u_{a_2})\geq E_{a_1}(\bar u)\geq e(a_1),$$
which  yields that  for all $a_2\searrow a_1$,
$$u_{a_2}\to \bar u\in \h\text{ and }\bar u\in \Lambda_{a_1}.$$
We thus obtain from (\ref{0.12}) that
\begin{equation}\label{0.13}-\frac{1}{2}\inte |\bar u|^4 dx\leq \liminf_{a_2\searrow a_1 }\frac{e(a_2)-e(a_1)}{a_2-a_1}\leq \limsup_{a_2\searrow a_1}\frac{e(a_2)-e(a_1)}{a_2-a_1}\leq  -\frac{1}{2}\gamma_{a_1}. \end{equation}
On the other hand,  by  (\ref{0.3}) we always have
$$\inte |\bar u|^4 dx\leq \gamma_{a_1}.$$
Thus, all inequalities in (\ref{0.13}) are indeed identities, from which we obtain
$$e'_+(a_1)=-\frac{1}{2}\gamma_{a_1}.$$
Similarly, if $a_2<a_1<a^*$,  letting $a_2\to a_1^-$ and repeating the above arguments, one can deduce that
$$e'_-(a_1)=-\frac{1}{2}\alpha_{a_1}.$$
This completes   the proof of Lemma \ref{lem2.B}.\qed
\begin{rem}\label{rem2.1}
Lemma \ref{lem2.B} implies that if $e(a)$ has a unique nonnegative minimizer, then $e(a)\in C^1([0,a^*))$. However, this is true generally for  $a\in[0,a_{*})$, where $a_{*}>0$ is given by (\ref{1.15}), in view of the possible multiplicity of nonnegative minimizers as $a$ approaches to $a^*$, cf. \cite{GS,GZZ}.
\end{rem}

Applying Lemma \ref{lem2.B} and Remark \ref{rem2.1}, we now prove Theorem \ref{pro1}.

\noindent\textbf{Proof of Theorem \ref{pro1}.} From  (\ref{0.8}), we have
\begin{equation*}
\begin{split}
|e(a_2)-e(a_1)|&\leq \frac{1}{2}|a_2-a_1|\max\Big\{\inte |u_{a_2}|^4dx,\inte |u_{a_1}|^4dx\Big\}\\
&\leq\frac{1}{2}\max\{\gamma_{a_1},\gamma_{a_2}\}|a_2-a_1|,\ \text{for all } a_1,a_2\in[0,a^*),
\end{split}
\end{equation*}
where $\gamma_{a_i}$ ($i=1,2$) is given by (\ref{0.3}).
 This implies that  $e(a)$ is locally  Lipschitz  continuous in $[0,a^*)$. It then follows  from Rademacher's theorem that $e(a)$ is differentiable for a.e. $a\in [0,a^*)$.
Moreover,  by  Lemma \ref{lem2.B} and Remark \ref{rem2.1} we see that
\begin{equation} \label{2.19c}
e'(a)\text{ exists for all $a\in[0,a_{*})$ and a.e. } a\in[0,a^*),\, \mbox{and} \ e'(a)=-\frac{1}{2}\inte |u|^4dx,\  \forall\ u \in\Lambda_a,
\end{equation}
and thus all  minimizers of $e(a)$ have the same $L^4(\R^2)$-norm.
Taking each nonnegative function $u_a\in \Lambda_a$, where $a\in[0,a^*)$ such that $e'(a)$ satisfies  (\ref{2.19c}), then $u_a $ satisfies  (\ref{d1.1}) for some Lagrange multiplier $\mu_a\in\R$. One can easily deduce from (\ref{d1.1})  and (\ref{2.19c}) that
\begin{equation}\label{2.18}
\mu_a=e(a)-\frac{a}{2}\inte |u_a|^4dx=e(a)+ae'(a).
\end{equation}
This shows that $\mu_a$ depends only on $a$ and is  independent of the choice of $u_a$. Thus for any given $a\in [0,a_*)$ and a.e. $a\in  [a_*,a^*)$, all minimizers of  $e(a)$ satisfy   equation (\ref{d1.1}) with the same Lagrange multiplier $\mu_a$.

We next prove the relationship $\Lambda_a=G_{a,\mu_a}$ to finish the proof of Theorem \ref{pro1}. For any given $a\in [0,a^*)$, consider any $u_a \in \Lambda_a$ and $u\in G_{a,\mu_a}$. Since  $u$ satisfies  (\ref{d1.1}) with $\mu=\mu_a$, we know that
  \begin{equation}\label{3.2.1}
\inte \Big[|\nabla u|^2+\big(V(x)-\mu_a\big)|u|^2\Big ]dx=a \inte |u|^4dx,
\end{equation}
which, together with (\ref{3.1.5}), implies that
 \begin{equation*}
J_{a,\mu_a}(u)=\frac{a}{4}\inte |u|^4dx.
\end{equation*}
Similarly,
\begin{equation}\label{3.2.2}J_{a,\mu_a}(u_a)=\frac{a}{4}\inte |u_a|^4dx.\end{equation}
Since $u\in G_{a,\mu_a}$, we have
 \begin{equation}\label{3.2.3}
J_{a,\mu_a}(u)=\frac{a}{4}\inte |u|^4dx\leq J_{a,\mu_a}(u_a)=\frac{a}{4}\inte |u_a|^4dx.
\end{equation}

Let $$\tilde u=\frac{1}{\sqrt\rho}u \text{ with}\ \rho=\inte |u|^2dx,$$
so that $\inte |\tilde u|^2 dx=1$. Noting that  $u_a$ is a minimizer of (\ref{eq1.3}), we obtain  that
$$E_a(\tilde u)\geq E_a(u_a).$$
Therefore, \begin{equation}\label{3.2.4}J_{a,\mu_a}(\tilde u)=\frac{1}{2}E_a(\tilde u)-\frac{\mu_a}{2}\inte|\tilde u|^2 dx\geq \frac{1}{2}E_a(u_a)-\frac{\mu_a}{2}\inte| u_a|^2 dx=J_{a,\mu_a}(u_a).\end{equation}
Note from (\ref{3.2.1}) that
\begin{align}
J_{a,\mu_a}(\tilde u)&=\frac{1}{2\rho}\inte \Big[|\nabla u|^2+\big(V(x)-\mu_a\big)|u|^2\Big]dx-\frac{a}{4\rho^2}\inte |u|^4dx\nonumber\\
&=\frac{a}{4\rho}\Big(2-\frac{1}{\rho}\Big)\inte |u|^4dx.\label{3.2.5}
\end{align}
It then follows from (\ref{3.2.2})-(\ref{3.2.5}) that
\begin{equation}\label{3.2.6}\frac{a}{4}\inte |u|^4dx=J_{a,\mu_a}(u)\leq J_{a,\mu_a}(u_a)\leq J_{a,\mu_a}(\tilde u)= \frac{a}{4\rho}\Big(2-\frac{1}{\rho}\Big)\inte |u|^4dx.\end{equation}
Since $a>0$ and $u\not=0$,  (\ref{3.2.6}) implies that
\begin{equation}\label{3.2.7}1\leq \frac{1}{\rho}\Big(2-\frac{1}{\rho}\big).\end{equation}
We thus obtain from  (\ref{3.2.7}) that  $\rho=1$, i.e., $\inte|u|^2dx=1$.  Moreover,   (\ref{3.2.6}) and  (\ref{3.2.7}) are identities,  i.e.,
$$J_{a,\mu_a}(u)=J_{a,\mu_a}(u_a),\ E_a(u)=E_a(u_a).$$
This implies that $u \in \Lambda_a$ and $u_a\in G_{a,\mu_a}$, and the proof is complete.\qed

We finally address the proof of Corollary \ref{pro1.4}, which deals with the uniqueness of nonnegative minimizers for $V(x)$ satisfying (\ref{1.17}).

\noindent{\textbf{Proof of  Corollary \ref{pro1.4}.}}
%Recall from Theorem 1.1 in \cite{GZZ} that $e(a)$ admits a unique nonnegative minimizer  $u_a$ for all $0<a<a_{*}$, where $a_{*}>0$ is given by (\ref{1.15}). Note that $u_a$ satisfies the Euler-Lagrange equation
%\begin{equation}\label{2.19}
%-\Delta u(x)+V(x)u(x)-\mu_{a} u(x)-a u^3(x)=0\quad  \text{ in }\
%\mathbb{R}^2
%\end{equation}
%for some suitable Lagrange multiplier $\mu _a\in \R$. By applying Theorem 2 in \cite{Li}, we derive from (\ref{2.19}) and  (\ref{1.17}) that $u_a$ must be radially symmetric about some point $x_0\in \R^2$ (especially, if $V$ is strictly increasing in $|x|$, then $u_a$  must be radially symmetric about the origin). In view of above results, to finish the proof of Theorem \ref{pro1.4} we next need only  to prove the results  for  a.e. $ a\in[a_{*},a^*)$ under the assumption (\ref{1.17}).
%
%Note from Theorem \ref{pro1} that for a.e. $a\in[a_{*},a^*)$,  all minimizers of $e(a)$ satisfy the same Euler-Lagrange equation (\ref{2.19}) with the same Lagrange multiplier $\mu _a\in \R$. Moreover, because   $V(x)$ satisfies (\ref{1.17}), it then follows again from Theorem 2 in \cite{Li} that $u_a$ solving  (\ref{2.19}) must be radially symmetric about $0\in \R^2$, and $u_a'(r)<0$ in $r=|x|>0$.  Further, applying Proposition 4.1  and Theorem 1.1 in \cite{BO} yields that positive radial solutions of (\ref{2.19}) must be unique. We thus conclude that for a.e. $a\in[a_{*},a^*)$, nonnegative minimizers of $e(a)$ must be unique and radially symmetric about  $0\in \R^2$, and the proof is therefore complete.\qed
Recall from Theorem 1.1 in \cite{GZZ} that $e(a)$ admits a unique nonnegative minimizer  $u_a$ for all $0<a<a_{*}$, where $a_{*}>0$ is given by (\ref{1.15}). Then by rotation we can deduce that  $u_a$ must be radially symmetric about the origin.  We next need only  to prove the results  for  a.e. $ a\in[a_{*},a^*)$ under the assumption (\ref{1.17}).

Note from Theorem \ref{pro1} that for a.e. $a\in[a_{*},a^*)$,  all minimizers of $e(a)$ satisfy the following Euler-Lagrange equation (\ref{2.19}) with the same Lagrange multiplier $\mu _a\in \R$,
\begin{equation}\label{2.19}
-\Delta u(x)+V(x)u(x)-\mu_{a} u(x)-a u^3(x)=0\quad  \text{ in }\
\mathbb{R}^2.
\end{equation}
Moreover, because   $V(x)$ satisfies (\ref{1.17}), it then follows again from Theorem 2 in \cite{Li} that $u_a$ solving  (\ref{2.19}) must be radially symmetric about $0\in \R^2$, and $u_a'(r)<0$ in $r=|x|>0$.  Further, applying Proposition 4.1  and Theorem 1.1 in \cite{BO} yields that positive radial solutions of (\ref{2.19}) must be unique. We thus conclude that for a.e. $a\in[a_{*},a^*)$, nonnegative minimizers of $e(a)$ must be unique and radially symmetric about  $0\in \R^2$, and the proof is therefore complete.\qed

%\begin{rem}
%From the proof of Theorem \ref{pro1.4}, we see that if $e'(a)$ exists, then all minimizers of  $e(a)$ satisfy the same Euler-Lagrange equation for all all $a\in[0,a_{**})$ and a.e. $a\in[a_{**},a^*)$, all  minimizers of
%\end{rem}

\section{General results on the concentration behavior for $u_a\in\Lambda_a$}\label{sec3}

In this section we prove Theorem \ref{thm1.2}, which focusses on the  concentration behavior of nonnegative minimizers for $e(a)$ as $a\nearrow a^*$ under general trapping potentials. Let $u_a$ be a non-negative minimizer of (\ref{eq1.3}). Then $u_a$ satisfies
the Euler-Lagrange equation
\begin{equation}\label{eq2.18}
-\Delta u_a(x)+V(x)u_a(x)=\mu_a u_a(x)+a u_a^3(x)\quad  \text{ in }\
\mathbb{R}^2,
\end{equation}
where $\mu_a\in \mathbb{R}$ is a suitable Lagrange multiplier. We first establish the following lemma.

\begin{lem}
Suppose $V(x)$ satisfies $(V_1)$.
%\begin{equation}\label{3.2}
%V(x)\in C(\R^2),\quad
%\underset{|x|\to\infty}{\lim} V(x) = \infty\, \text{ and }\, \underset{x\in \R^2}{\inf}
%V(x) =0.
%\end{equation}
Let $u_a\in\Lambda_a$ and let $\varepsilon _a>0$ be defined as
\begin{equation}\label{eq2.9a}\varepsilon ^{-2}_a:=\int_{\mathbb{R}^2}|\nabla u_a(x)|^2dx.
\end{equation}
Then,
\begin{itemize}
  \item [\rm(i)] \begin{equation}\label{3:nab}
\eps_a\to 0\quad \text{as}\quad a\nearrow a^*\,.
 \end{equation}
  \item [\rm(ii)]  $u_a(x)$ has at least one local maximum point $\bar z_a$, and there exists $\eta>0$ such that the
normalized function
\begin{equation}\label{3.5}\bar{w}_a(x)=\varepsilon_au_a(\varepsilon_a x+\bar{z}_a)\end{equation}
satisfies
 \begin{equation}\label{eq2.12}
\int_{\mathbb{R}^2}|\nabla \bar w_a|^2dx=1,\quad
\int_{\mathbb{R}^2}| \bar w_a|^4dx\to\frac{2}{a^*}\quad \text{as}\quad a\nearrow a^*,
\end{equation}
and
\begin{equation}\label{eq2.16}
\liminf_{a\nearrow a^*}\int_{B_{2}(0)}|\bar w_a|^2dx\geq\eta>0.
\end{equation}
\item [\rm(iii)] The local maximum point $\bar z_a$ of $u_a(x)$ satisfies
\begin{equation}\label{new3.7}
\lim_{a\nearrow a^*} {\rm dist} (\bar z_a, \mathcal{A})=0,
\end{equation}
where $\mathcal{A}=\{x\in\R^2: V(x)=0\}$ denotes the zero point set of $V(x)$.
\end{itemize}
\end{lem}

\noindent \textbf{Proof.} \textbf{(i).}
If (\ref{3:nab}) is false, then there exists a sequence $\{a_k\}$, where $a_k\nearrow a^*$ as $k\to\infty$, such that  $\{u_{a_k}(x)\}$ is bounded uniformly in $\h$. By applying  the compactness of Lemma \ref{2:lem1},  there exist a subsequence (still denoted by $\{a_k\}$) of $\{a_k\}$ and $u_0\in \h$ such that
$$u_{a_k}\overset{k}\rightharpoonup u_0 \quad \text{weakly in }\h \quad \text{and}\quad u_{a_k}\overset{k}\rightarrow u_0\quad \text{strongly in }L^2(\R^2)\,.$$
Consequently,
  $$0=e(a^*)\leq E_{a^*}(u_0)\leq\lim_{k\to\infty} E_{a_k}(u_{a_k})=\lim_{k\to\infty} e(a_k)=0,$$
since Theorem 1 of \cite{GS}, $e(a)\to0$ as $a\nearrow a^*$. This then indicates that $u_0$ is a minimizer of $e(a^*)$, which is impossible since Theorem 1 of \cite{GS} shows that $e(a^*)$ cannot be attained. So, part \textbf{(i)} is proved.

\textbf{(ii).}
Since $u_a\geq0$ satisfies the equation (\ref{eq2.18}) and  $\underset{|x|\to\infty}\lim V(x)=+\infty$, for any fixed $a\in[0,a^*)$, one can use comparison principle as in \cite{KW} to deduce  that $u_a$ decays exponentially to zero  at infinity, thus $$u_a(x)\to0\quad\text{as}\quad |x|\to\infty.$$
This implies that each $u_a$ has at least one local maximum point, which is denoted by $\bar z_a$. Let $\bar w_a$ be defined by (\ref{3.5}),  we have
\begin{equation*}
0\leq\int_{\mathbb{R}^2}|\nabla
u_a(x)|^2dx-\frac{a}{2}\int_{\mathbb{R}^2}| u_a(x)|^4dx=\eps
^{-2}_a-\frac{a}{2}\int_{\mathbb{R}^2}| u_a(x)|^4dx\leq
e(a)\rightarrow 0\ \text{ as }\ a\nearrow a^*.
\end{equation*}
We thus obtain from (\ref{3:nab}) that
\begin{equation}\label{eq2.10}
\frac{a}{2}\int_{\mathbb{R}^2}|
u_a(x)|^4dx\cdot\eps ^{2}_a\to1\quad \text{ as }\quad
a\nearrow a^*.
\end{equation}
Together with  (\ref{eq2.9a}) and (\ref{3.5}), we conclude from the above that (\ref{eq2.12}) holds.

To complete the proof of (ii), it now remains to prove (\ref{eq2.16}). We first claim that
\begin{equation}\label{3.10}
 \bar w_a^2(0)\geq \frac{1}{2a^*}\quad \text{as}\quad a\nearrow a^*.
\end{equation}
Indeed, by (\ref{eq2.18}) we see   that
\begin{equation*}
\mu_a=e(a)-\frac{a}{2}\int_{\mathbb{R}^2}|u_a|^4dx.
\end{equation*}
It then follows from  (\ref{3:nab}) and (\ref{eq2.10}) that
\begin{equation}\label{3.12}\eps_a^2\mu_a\to-1\quad \text{as}\quad a\nearrow a^*.\end{equation}
Moreover, in view of (\ref{eq2.18}), $\bar w_a(x)$
satisfies the elliptic equation
\begin{equation}\label{eq2.19}
-\Delta \bar w_a(x)+\eps_a^2V(\eps_a x+\bar z_a
)\bar w_a(x)=\eps_a^2\mu_a \bar w_a(x)+a \bar w_a^3(x)\quad \text { in
} \ \mathbb{R}^2.
\end{equation}
Since $\bar w_a(x)$ attains its local maximum at $x=0$ and note that $V(x)\geq0$, we thus obtain from (\ref{eq2.19}) and (\ref{3.12}) that
$$\bar w_a^2(0)\geq \frac{-\eps_a^2\mu_a}{a}\geq \frac{1}{2a^*}\quad \text{as}\quad a\nearrow a^*,$$
i.e., (\ref{3.10}) holds and the claim is proved.

Note from (\ref{3.12}) and  (\ref{eq2.19}) that for $a\nearrow a^*$,
\begin{equation}\label{3.16}
-\Delta \bar w_a-c(x)\bar w_a\leq 0 \  \text{ in }\ \R^2,\  \text{ where }\ c(x)=a \bar w_a^2(x).\end{equation}
 By applying De Giorgi-Nash-Moser theory to (\ref{3.16}) (see \cite[Theorem 4.1]{HL}), we then
have  \begin{equation}\label{new3.15}\max_{B_1(\xi)} \bar w_a\leq
C\Big(\int_{B_2(\xi)}|\bar w_a|^2 dx\Big)^\frac{1}{2},\end{equation}
where
$\xi$ is an arbitrary point in $\mathbb{R}^2$, and $C>0$ depends only on the upper  bound of $\|c(x)\|_{L^2(B_2(\xi))}$, i.e.,  the upper  bound of $\|\bar w_a\|_{L^4(B_2(\xi))}$. Therefore, it then follows from (\ref{eq2.12}) that $C>0$ is bounded uniformly  as $a\nearrow a^*$.
Taking $\xi=0$, we thus obtain from (\ref{3.10}) and (\ref{new3.15}) that
$$\liminf_{a\nearrow a^*}\int_{B_2(0)}|\bar w_a|^2 dx\geq \frac{1}{2C^2a^*}:=\eta>0,$$
and  (\ref{eq2.16}) is therefore established.

\textbf{(iii).} By the Gagliardo-Nirenberg inequality
(\ref{GNineq}) and Theorem 1 of \cite{GS}, we have
\begin{equation}\label{eq2.8}
\int_{\mathbb{R}^2}V(x)|u_a(x)|^2dx\leq e(a)\to 0
\quad \text{ as }\quad a\nearrow a^*\,.
\end{equation} It then follows from (\ref{3.5}) and  (\ref{eq2.8}) that
\begin{equation}\label{eq2.17}
\int_{\mathbb{R}^2}V(x)|u_a(x)|^2dx=\int_{\mathbb{R}^2}V(\eps_a
x+\bar z_a)|\bar w_a(x)|^2dx\rightarrow 0 \quad \text{as}\quad
a\nearrow a^*.
\end{equation}
On the contrary, suppose (\ref{new3.7}) is false. Then there exist a constant $\delta >0$ and a subsequence $\{a_n\}$ with
$a_n\nearrow a^*$ as $n\to\infty $ such that
$$\eps_n:=\eps_{a_n} \overset{n}\to 0 \quad \text{and}\quad \lim_{n\to\infty}{\rm dist}(\bar z_{a_n}, \mathcal{A})\geq \delta>0.$$ Since $\underset{|x|\to\infty}\lim{V(x)}=\infty$, we thus deduce that there exists $C_\delta>0$ such that
 $$\lim_{n\to\infty}V(\bar z_{a_n})\geq2C_\delta>0.$$
Then, by Fatou's Lemma  and (\ref{eq2.16}), we derive that
\begin{equation*}
\lim_{n\rightarrow\infty}\int_{\mathbb{R}^2}V(\eps_n
x+\bar z_{a_n})|\bar w_{a_n}(x)|^2dx\geq
\int_{\mathbb{R}^2}\lim_{n\rightarrow\infty}V(\eps_n
x+\bar z_{a_n})|\bar w_{a_n}(x)|^2dx\geq C_\delta \eta,
\end{equation*}
which however contradicts (\ref{eq2.17}). Therefore, (\ref{new3.7}) holds and the proof is complete.\qed
\\

\noindent\textbf{Proof of Theorem \ref{thm1.2}.}
For any given sequence $\{a_k\}$ with  $a_k\nearrow a^*$ as
$k\to\infty$,   we denote $$u_k(x):=u_{a_k}(x),
\bar w_k:=\bar w_{a_k}\ge 0, \text{ and } \eps_k:=\eps_{a_k}>0,$$
 where $\eps_{a_k}$ is defined by (\ref{eq2.9a}) and satisfies $\eps_{a_k}\to  0$ as $k\to\infty$. Let $\bar z_k:=\bar z_{a_k}$  be a local maximum point of $u_k(x)$.  It  yields from (\ref{new3.7}) and (\ref{eq2.12}) that there exists a subsequence, still denoted by $\{a_k\}$,  of $\{a_k\}$ such that
\begin{equation}\label{3.23}\lim_{k\to\infty}\bar z_k= x_0\ \text{and }\ V(x_0)=0,
\end{equation}
and
 $$\bar w_k\overset{k}\rightharpoonup w_0\geq 0 \text{ weakly in } H^1(\mathbb{R}^2)$$
 for some $w_0\in H^1(\R^2)$. Moreover, $\bar w_k$ satisfies
\begin{equation}\label{3:222}
-\Delta \bar{w}_k(x)+\varepsilon_k^2 V(\varepsilon_k
x+\bar{z}_k)\bar{w}_k(x)=\mu_k\varepsilon_k^2 \bar{w}_k(x)+a _k
\bar{w}_k^{3}(x) \ \text { in } \ \mathbb{R}^2,
\end{equation}
and note from (\ref{3.12}) that $\mu_k\varepsilon_k^2 \to -1$ as $k\to\infty $. Thus, by
taking the  weak limit in (\ref{3:222}), we obtain that    $\bar w_0$ satisfies
\begin{equation}\label{eq2.20}
-\Delta w(x)=- w(x)+a^* w^{3}(x)\quad \text { \ in }
\ \mathbb{R}^2.
\end{equation}
Furthermore, it follows from (\ref{eq2.16}) that $\bar w_0\not\equiv 0$,
and thus  $\bar w_0>0$ by the strong maximum principle. By a
simple rescaling,  the uniqueness (up to translations) of positive
solutions for the nonlinear scalar field equation (\ref{Kwong})
implies that
\begin{equation}\label{eq2.21}
\bar w_0(x)=\frac{1}{\|Q\|_2}Q(|x-y_0|) \quad \text{for
some} \quad y_0\in\mathbb{R}^2,
\end{equation}
where $||\bar w_0||_2^2=1$. By the norm preservation we further conclude
that
$$\bar w_k \overset{k}{\to} \bar w_0 \ \text{ strongly in }L^2(\R^2).$$
Together with the boundness of $\bar w_k$ in $H^1(\mathbb{R}^2) $, this implies that
$$\bar w_k \overset k\to \bar w_0 \
 \text{ strongly in $L^p(\R^2)$ for any $2\leq p<\infty$}.$$
Moreover,  since $\bar w_k$ and
$\bar w_0$ satisfy (\ref{3:222}) and (\ref{eq2.20}), respectively, a
simple analysis shows that
\begin{equation}\label{3.15}
\bar w_k \overset{k}\to \bar w_0=\frac{1}{\|Q\|_2}Q(|x-y_0|)\, \text{ strongly in }\, H^1(\mathbb{R}^2).
\end{equation}

Since $V(x)\in C^\alpha(\R^2)$, one can deduce from (\ref{3:222}) that $\bar{w}_k\in C^{2,\alpha_1}_{loc}(\mathbb{R}^2)$ for some $\alpha_1\in(0,1)$, which yields that
\begin{equation}\label{3.22}\bar{w}_k\xrightarrow{k}
\bar{w}_0 \text{ in } C^2_{loc}(\mathbb{R}^2).\end{equation}
Since $x=0$
is a critical (local maximum) point of $\bar{w}_k(x)$ for all $k>0$,  in view of (\ref{3.22}) it is also a
critical point of $\bar{w}_0$.  We therefore conclude from the
uniqueness (up to translations) of positive radial solutions for
(\ref{Kwong}) that   $\bar{w}_0$ is spherically symmetric about the
origin, i.e. $y_0=(0,0)$ in (\ref{3.15}) and
\begin{equation}\label{eq2.350}
\bar w_0(x)=\frac{1}{\|Q\|_2}Q(|x|)\,.
\end{equation}

One can deduce from (\ref{3:222}) that $\bar{w}_k\geq
(\frac{1}{2a^*})^\frac{1}{2}$ at each local maximum point.
Since $\bar{w}_k$ decays to zero uniformly in $k$ as
$|x|\rightarrow\infty$, all   local maximum points of $\bar{w}_k$
stay in a finite ball in $\mathbb{R}^2$. It then follows
from (\ref{3.22})  and Lemma 4.2 in \cite{NT} that for large $k$, $\bar{w}_k$ has no
critical points other than the origin. This  gives the uniqueness of
local maximum points
for $\bar{w}_k(x)$, which therefore implies that $\bar z_k$ is the  unique maximum point of $u_k$  and $\bar z_k$ goes to a  global minimum point of $V(x)$ as $k\to\infty$. Moreover, (\ref{con:a1}) is followed from
(\ref{3.15}) and (\ref{eq2.350}).  The proof of Theorem \ref{thm1.2} is therefore complete.
\qed

\section{ Concentration behavior under multiple-well  potential}\label{sec4}
In this section we turn to proving Theorem \ref{thm1.3}, which gives more detailed concentration behavior for the  nonnegative minimizers of $e(a)$ as $a\nearrow a^*$, provided that the  potential satisfies $(V_1)$ and $(V_2)$, that is, multiple-well potential.
% $R$ and $\Lambda$ are defined as
%\begin{equation}\label{4.1}
%R=\max_{i=1,\cdots, m}R_i,\ \text{ where }\  R_i:=\max_{x\in\Omega_i}dist (x,\partial \Omega_i)>0\,,
%\end{equation}
%and \begin{equation}\label{4.2}
%\Lambda:=\big\{\Omega_i:\, R_i=R\ \mbox{for some} \ i\in 1,\cdots,m\big\}\,.
%\end{equation}
Inspired  by \cite{GS,GZZ}, we start with establishing the following energy estimates of $e(a)$.

\begin{lem}\label{lem2.2}
Suppose $V(x)$ satisfies  $(V_1)$ and $(V_2)$, then
\begin{equation}\label{2:en}
0\leq e(a)\leq \frac{1+o(1)}{4R^2a^*}(a^*-a)\big(|\ln(a^*-a)|\big)^2 \text{ as}\quad  a\nearrow a^*\,.
\end{equation}
\end{lem}

\noindent{\bf Proof.}  Let $\varphi(x)\in
C^\infty_0(\R^2)$ be
 a nonnegative smooth cut-off function such that $\varphi(x)=1$ if $|x|\leq1$
 and $\varphi(x)=0$ if $|x|\geq2$. Choose $\Omega_i\in \Lambda$ and  $x_0\in\Omega_i$ for some $i\in\{1,\cdots ,m\}$ such that
  \begin{equation}\label{2:V}{\rm dist} (x_0, \partial \Omega_i)=R_i=R\,,\end{equation}
 where  $R$ and $\Lambda$ are defined by (\ref{1:R}) and (\ref{1.18}), respectively.

 Let
 \begin{equation}\label{2:trial}
\phi(x)=A_{\tau R}\frac{\tau}{\|Q\|_2}\varphi\big(\frac{x-
x_0}{R}\big)Q\big(\tau|x-x_0|\big)\text{ for }\tau>0\,,
 \end{equation}
where $A_{\tau R}>0$  is chosen such that $\int_{\R^2}\phi^2dx=1$.
By the exponential decay of $Q$ in (\ref{4:exp}), we obtain that
\begin{equation}\label{2:limt1}
1\le A_{\tau R}^2\leq 1+Ce^{-2\tau R}
\quad\text{as}\quad \tau\to\infty\,,
\end{equation}
\begin{eqnarray}\label{2:nabla}
&&\int_{\R^2}|\nabla \phi|^2dx\leq\frac{A_{\tau R}^2\tau^2}{{a^*}}\int_{\R^2}|\nabla Q|^2dx+C\tau e^{-2\tau R}\quad\text{as}\quad \tau\to\infty\,,
\end{eqnarray}
and
\begin{equation}\label{2:L4}
\inte|\phi|^4dx\geq \frac{A_{\tau R}^4\tau^2}{{a^*}^2}\inte|Q|^4dx-C e^{-3\tau R}\quad\text{as}\quad \tau\to\infty\,.
\end{equation}
Using the equality (\ref{1:id}), we then derive from (\ref{2:limt1})-(\ref{2:L4}) that
\begin{equation}\label{2:limt2}
 \begin{split}
 &\int_{\R^2}|\nabla \phi|^2dx-\frac{a}{2}\int_{\R^2} |\phi|^4dx\\
&\leq \frac{A_{\tau R}^2\tau^2}{a^*}\int_{\R^2}|\nabla
Q|^2dx-\frac{A_{\tau R}^2a\tau^2}{2{a^*}^2}\int_{\R^2} Q^4dx+C\tau e^{-2\tau R} \\
&= A_{\tau R}^2\frac{a^*-a}{a^*}\tau^2+C\tau e^{-2\tau R}\quad\text{as}\quad \tau\to\infty\,.
\end{split}\end{equation}
On the other hand, since
$$V\big(\frac{x}{\tau}+x_0\big)\leq \sup_{y\in B_{2R}(x_0)}V(y) \quad \text{ for all }\tau>0 \text{ and }x\in B_{2\tau R}(0)\,,$$
we deduce from (\ref{2:V})  that
$$V\big(\frac{x}{\tau}+x_0\big)\equiv0\quad \text{ for all } x\in B_{\tau R}(0)\,.$$
It then follows  from the exponential decay of $Q$ that
\begin{equation*}
\begin{split}
\inte V(x)\phi^2(x)dx&=\frac{A_{\tau R}^2}{\|Q\|_2^2}\int_{B_{2\tau R}(0)}V\big(\frac{x}{\tau}+x_0\big)
\varphi^2\big(\frac{x}{\tau R}\big)Q^2(x)dx\\
&\leq C\int_{B_{2\tau R}(0)\setminus B_{\tau R}(0)}|Q(x)|^2dx\leq Ce^{-2\tau R}\quad\text{as}\quad \tau\to\infty\,.
\end{split}
\end{equation*}
Combining with (\ref{2:limt2}), this implies that
\begin{equation}
0\le E_a(\phi)\leq A_{\tau R}^2\frac{a^*-a}{a^*}\tau^2+C\tau e^{-2\tau R}\quad\text{as}\quad \tau\to\infty\,.
\end{equation}
By setting $\tau=\frac{1}{2R}|\ln(a^*-a)|$, we then conclude from the above that
\begin{equation}
\begin{split}
0\le E_a(\phi)&\leq  \frac{A_{\tau R}^2}{4R^2a^*}(a^*-a)\big(|\ln(a^*-a)|\big)^2+\frac{C}{2R}(a^*-a)|\ln(a^*-a)|\\
&=\frac{1+o(1)}{4R^2a^*}(a^*-a)\big(|\ln(a^*-a)|\big)^2 \quad \mbox{as}\quad a\nearrow a^*,
\end{split}
\end{equation}
and the proof is therefore done.
\qed

In order to derive the optimal energy estimates of $e(a)$ as $a\nearrow a^*$, we also need the following lemma.

\begin{lem}\label{lem2.3}
Under the assumption   $(V_1)$, suppose $w_\eps(x)\in \h$ is a nonnegative  solution of the following equation
\begin{equation}\label{2:eps}
-\Delta w_{\eps}+\eps^2V(\eps x+x_\eps)w_\eps=-\beta^2_\eps w_\eps+a_\eps w_\eps^3\quad \text{in}\quad \R^2,\ \eps>0\,,
\end{equation}
where both  $x_\eps\in\R^2$ and $a_\eps\geq0$ are bounded uniformly as $\eps\to 0$, and $w_\eps$ satisfies \begin{equation}\label{2:h}w_\eps(x) \rightarrow w(x)>0\quad \text{strongly in}\quad H^1(\R^2)\quad \text{and}\quad \beta_\eps\to \beta>0\quad \text{as}\quad\eps\to 0\,.\end{equation}
Then for any $\sigma\in(0,1)$ and $\bar R>0$, there exist $\varepsilon_0=\varepsilon_0(\beta,\sigma,\bar R)>0$ and $\bar\mu=\bar\mu(\beta,\sigma)>0$ such that for all $\eps\in(0,\eps_0)$,
\begin{equation}
w_\eps(x)\geq \bar\mu e^{-(\beta+\sigma)|x|} \quad \text{for}\quad \rho(\beta,\sigma)\leq|x|\leq \frac{\bar R}{\eps}\,,
\end{equation}
where \begin{equation}\label{2:rho}
\rho(\beta,\sigma):=\max\Big\{(\beta+\sigma)^{-1}\big(\frac{3-\sqrt5}{2}\big)^{-2},\frac{\beta+\sigma}
{\sigma(\beta+\frac34\sigma)}\Big\}\,.\end{equation}
\end{lem}

\noindent{\bf Proof.}
For any $\sigma>0$ and $\bar R>0$, denote $R_\eps:=\frac{2\bar R}{\eps}>0$.  Since $V(x)$ satisfies $(V_1)$ and  $\{x_\eps\}$ is bounded uniformly as $\eps\to0$, we have
$$\eps^2V(\eps x+x_\eps)\leq C\eps^2\quad \text{for}\quad x\in B_{R_\eps+4}(0)\,.$$
In addition, since  $w_\eps\geq0$ in $\R^2$ and $\beta_\eps\to\beta>0$ as $\eps\to0$, it then follows from (\ref{2:eps}) that there exists $\eps_1=\eps_1(\beta, \sigma, \bar R)>0$ such that for all $\eps\in(0,\eps_1)$,
\begin{equation}\label{2:super}
-\Delta w_\eps(x)+\big(\beta+\frac{\sigma}{2}\big)^2w_\eps(x)\geq0 \quad\text{in}\quad  B_{R_\eps+4}(0) \,.
\end{equation}
By applying  Theorem 8.18 in  \cite{GT}  and (\ref{2:h}), we derive from the above that for any $p>1$ and $\eps\in(0,\eps_1)$,
\begin{equation}\label{eq4.15}\inf_{B_1(\xi)}w_\eps(x)\geq C\|w_\eps\|_{L^p(B_2(\xi))}\geq \frac{C}{2}\|w\|_{L^p(B_2(\xi))}\quad \text{for all}\quad \xi\in B_{R_\eps}(0)\,.\end{equation}
Moreover, we note from (\ref{2:rho}) that there exists $\eps_2=\eps_2(\beta, \sigma, \bar R)>0$ such that $\rho(\beta,\sigma)<\frac{R_\eps}{2}$ for all $\eps\in(0,\eps_2)$,  and thus by letting
\begin{equation}\label{eq4.17}\eps_0=\eps_0(\beta,\sigma, \bar R):=\min\{\eps_1(\beta,\sigma, \bar R),\eps_2(\beta,\sigma, \bar R)\},\end{equation}
it then follows from  (\ref{eq4.15}) that  there exists  a constant  $C(\beta,\sigma)>0$ such that
\begin{equation}\label{eq4.16}w_\eps(x)\geq  \frac{C}{2}\|w\|_{L^p(B_2(x))}\geq C(\beta,\sigma)>0\quad \text{for all } |x|= \rho(\beta,\sigma) \text{ and }\eps\in(0,\eps_0).\end{equation}

Let
\begin{equation}\label{2:test}
U(x)=\mu e^{-\lambda|x|},\end{equation}
where $0<\mu\leq C(\beta,\sigma)$ and $\lambda>0$ to be determined later.
It is easy to check that
\begin{equation}\label{2:test1}
-\Delta U(x)+\lambda^2U(x)=\frac{\lambda U(x)}{|x|}\quad \mbox{in}\quad \R^2\,.
\end{equation}
Let $\varphi(x)\in C_0^\infty(\R^2)$ be given by
\begin{equation}\label{2:smo}\varphi(x)=\begin{cases}
\exp\{\frac{1}{|x|^2-1}\}\quad &|x|<1,\\
0\quad &|x|\geq1\,,
\end{cases}\end{equation}
and define $$h_\eps(x):=\varphi_{R_\eps}(x) U(x)=\varphi(\frac{x}{R_\eps})U(x)\,.$$
It follows from  (\ref{2:test}) and (\ref{2:test1})  that, for $x'=\frac{x}{R_\eps}$,
\begin{equation*}
\begin{split}
\Delta h_\eps&=\varphi_{R_\eps}(x)\Delta U+\frac{2}{R_\eps}\nabla U \cdot\nabla _{x'}\varphi(\frac{x}{R_\eps})+\frac{U}{R_\eps^2}\Delta_{x'}\varphi(\frac{x}{R_\eps}) \\
&=\Big(\lambda^2U(x)-\frac{\lambda U(x)}{|x|}\Big)\varphi_{R_\eps}(x)-\frac{2\lambda Ux\cdot \nabla_{x'}\varphi(\frac{x}{R_\eps}) }{R_\eps |x|}+\frac{U}{R_\eps^2}\Delta_{x'}\varphi(\frac{x}{R_\eps})\quad \mbox{in}\quad \R^2,
\end{split}
\end{equation*}
i.e.,
\begin{equation}\label{5.2.18}
-\Delta h_\eps+\lambda^2h_\eps=\frac{\lambda h_\eps}{|x|}+\frac{2\lambda Ux\cdot \nabla_{x'}\varphi(\frac{x}{R_\eps}) }{R_\eps |x|}-\frac{U}{R_\eps^2}\Delta_{x'}\varphi(\frac{x}{R_\eps})\quad \mbox{in}\quad \R^2.
\end{equation}
By (\ref{2:smo}), we have, for $i=1,2$,
\begin{equation}\label{5.2.19}D_i\varphi(x)=\begin{cases}
-\varphi(x)\frac{2x_i}{(|x|^2-1)^2}\quad &|x|<1,\\
0\quad &|x|\geq1\,.
\end{cases}\end{equation}
Direct calculations show that
\begin{equation}\label{5.2.20}\Delta\varphi(x)=4\varphi(x)\frac{|x|^4+|x|^2-1}{(|x|^2-1)^4}\quad \mbox{in}\quad \R^2.\end{equation}
It then follows from (\ref{5.2.19}) and (\ref{5.2.20}) that
 \begin{equation}\label{5.2.21}\frac{2\lambda Ux\cdot \nabla_{x'}\varphi(\frac{x}{R_\eps}) }{R_\eps |x|}=-\frac{4\lambda h_\eps(x)R_\eps^2|x|}{(|x|^2-R_\eps^2)^2},\end{equation}
 and
 \begin{equation}\label{5.2.22}
 \frac{U}{R_\eps^2}\Delta_{x'}\varphi(\frac{x}{R_\eps})
 =\frac{4h_\eps(x)R_\eps^2[|x|^4+|x|^2R_\eps ^2-R_\eps ^4]}{(|x|^2-R_\eps ^2)^4}.\end{equation}
Putting  (\ref{5.2.21}) and (\ref{5.2.22}) into (\ref{5.2.18}), we then have
\begin{equation}\label{2:test2}
-\Delta h_\eps(x)+\lambda^2h_\eps(x)=\frac{\lambda h_\eps(x)}{|x|}-\frac{4h_\eps(x)R_\eps^2}{(|x|^2-R_\eps^2)^4}\Big(\lambda |x| (|x|^2-R_\eps^2)^2+|x|^4+R_\eps^2|x|^2-R_\eps^4\Big)
\end{equation}
in $\R^2$. Note that
$$|x|^4+R_\eps^2|x|^2-R_\eps^4\geq0,\ \text{ if }\ \frac{-1+\sqrt5}{2}R_\eps^2\leq|x|^2\leq R_\eps^2\,,$$
and
$$\lambda|x|(|x|^2-R_\eps^2)^2\geq R_\eps^4,\quad \text{if }\ \lambda^{-2}\big(\frac{3-\sqrt5}{2}\big)^{-4}\leq|x|^2<\frac{-1+\sqrt5}{2}R_\eps^2\,.$$
Therefore,
$$\lambda|x| (|x|^2-R_\eps^2)^2+|x|^4+R_\eps^2|x|^2-R_\eps^4\geq0,\quad \text{if}\quad \lambda^{-1}\big(\frac{3-\sqrt5}{2}\big)^{-2}\leq|x|\leq R_\eps\,.$$
This estimate and (\ref{2:test2}) yield that
\begin{equation*}
-\Delta h_\eps(x)+\lambda^2h_\eps(x)\leq\frac{\lambda h_\eps(x)}{|x|}
\quad \text{for}\quad \lambda^{-1}\big(\frac{3-\sqrt5}{2}\big)^{-2}\leq|x|<R_\eps\,.
\end{equation*}
Taking $\lambda=\beta+\sigma$, we then further derive from (\ref{2:rho}) that for all $\eps\in(0,\eps_0)$ with $\eps_0$ given by (\ref{eq4.17}),
\begin{equation}\label{2:sub}
\begin{split}
& -\Delta h_\eps(x)+\big(\beta+\frac{\sigma}{2}\big)^2h_\eps(x)\\
\leq &\Big[\frac{\lambda}{|x|}-\sigma\big(\beta+\frac34\sigma\big)\Big]h_\eps(x)\leq0\quad \text{for}\quad \rho(\beta,\sigma)\leq|x|<R_\eps\,.
\end{split}\end{equation}

Since   $0<\mu<C(\beta,\sigma)$, where $\mu$ is as in (\ref{2:test}),  it then  follows from (\ref{eq4.16}) that  $w_\eps(x)\geq h_\eps(x)$ at
$ |x|=\rho(\beta,\sigma)$. Moreover,  $h_\eps(x)\equiv0\leq  w_\eps(x)$ at $|x|=R_\eps$. We thus obtain from (\ref{2:super}) and (\ref{2:sub}) that for all $\eps\in(0,\eps_0)$,
$$w_\eps(x)\geq h_\eps(x)=\varphi_{R_\eps}(x)\mu e^{-(\beta+\sigma)|x|}\quad \text{for}\quad \rho(\beta,\sigma)\leq|x|<R_\eps\,.$$
Since  $\varphi_{R_\eps}(x)\geq e^{-\frac{4}{3}}$ for $|x|\leq\frac{R_\eps}{2}=\frac{\bar R}{\eps}$, we conclude that
$$w_\eps(x)\geq\mu e^{-\frac{4}{3}} e^{-(\beta+\sigma)|x|}\quad \text{for}\quad \rho(\beta,\sigma)\leq|x|\leq\frac{\bar R}\eps\,,$$
and the proof is done.
\qed

Following above lemmas and Theorem \ref{thm1.2}, we next complete the proof of Theorem \ref{thm1.3}.

\noindent\textbf{Proof of Theorem \ref{thm1.3}.}
For any  sequence
$\{a_k\}$ with $a_k\nearrow a^*$ as $k\to\infty$, we still denote $u_k:=u_{a_k}$  a nonnegative minimizer of $e(a_k)$.  It then follows from Theorem \ref{thm1.2} that
there exists a subsequence of $\{a_k\}$, still denoted by $\{a_k\}$,  such that each $u_k$ has a
unique maximum point $\bar z_k$, which satisfies $\bar z_k\overset{k}{\to} x_0$ for some $x_0\in\bar\Omega$.
Moreover, we have
\begin{equation*}%\label{1:na1}
\varepsilon_k:=\Big(\inte|\nabla u_{a_k}|^2dx\Big)^{-\frac{1}{2}}\to0\quad \text{ as }\ k\to\infty\,,\end{equation*}
 and
\begin{equation}\label{4.28}
\bar w_k(x):=\varepsilon_k u_{k}\big(\varepsilon_k x+\bar z_k
\big) \to \frac { Q(x)}{\sqrt{a^*}}\ \text{ strongly in }\ H^1(\R^2)
\end{equation}
as $k\to\infty$. To complete the proof of Theorem \ref{thm1.3}, it thus remains to establish  the estimates (\ref{5.1.4})-(\ref{con:a}).

Since $x_0\in\bar\Omega$,  without loss of generality we may assume that $x_0\in\bar\Omega_{i_0}$ for some $i_0\in\{1,\ldots,m\}$. We first claim that (\ref{1:R0}) holds.
Indeed, since $\bar z_k\to x_0$ as $k\to\infty$, where $x_0\in\bar\Omega_{i_0}$, we  always have
 \begin{equation}\label{3:rel}
 R_0:=\lim_{k\to\infty}{\rm dist}(\bar z_k, \partial \Omega_{i_0})\leq R_{i_0}\leq R\,.\end{equation}
For any fixed $0<\delta<\min\{\frac{R_{i_0}}{5},\frac{d}{5}\}$, where $d>0$ is defined by
\begin{equation*}%\label{1:dist}
d=\inf\big\{{\rm dist}(\Omega_i,\Omega_j);\, i\neq j\,\text{ and } \, i,j=1,\cdots, m\big\}\,.
\end{equation*}
Since $\bar z_k\to x_0\in\bar \Omega_{i_0}$ as $k\to\infty$,  we obtain that
$\bar z_k\in \Omega_{i_0}^\delta$ for large $k$ and
\begin{equation}\label{3:dist2}{\rm dist}(x_0, \partial \Omega_{i_0})=R_0\,,\end{equation}
where $\Omega_{i_0}^\delta:=\{x\in\R^2: {\rm dist}(x,\Omega_{i_0})<\delta\}$.
Choosing $$r_1=R_0+3\delta, \quad r_2=R_0+4\delta\,,$$
 it then follows from (\ref{3:dist2}) that
$$\Big|\big(B_{r_2}(x_0)\setminus B_{r_1}(x_0)\big)\cap (\Omega^\delta)^c\Big|>0\,.$$
Recalling that  $\bar z_k\to x_0$ as $k\to\infty$, then for large $k>0$,
\begin{equation}\label{3:bound}
\Big|\big(B_{r_2}(\bar z_k)\setminus B_{r_1}(\bar z_k)\big)\cap (\Omega^\delta)^c\Big|\geq\frac{1}{2}\Big|\big(B_{r_2}(x_0)\setminus B_{r_1}(x_0)\big)\cap (\Omega^\delta)^c\Big|>0\,.\end{equation}
In view of (\ref{3.12})  and  $\bar w_k(x)$ satisfies (\ref{3:222}) and (\ref{4.28}),   we can apply  Lemma \ref{lem2.3} with $\beta=1$ to see that for any $\sigma\in(0,1)$,
\begin{equation}\label{3:bound1}\bar w_k(x)\geq C(\sigma) e^{-(1+\sigma)|x|}\quad \text{in }\ B_{\frac{r_2}{\eps_k}}(0)\setminus B_{\frac{r_1}{\eps_k}}(0)\end{equation}
holds for large $k>0$.
Moreover, note that there exists $C=C(\delta)>0$ such that
\begin{equation}\label{3:bound2}
V(x)\geq C\quad \text{for}\quad x\in\big(B_{r_2}(\bar z_k)\setminus B_{r_1}(\bar z_k)\big)\cap (\Omega^\delta)^c\,.
\end{equation}
We thus derive from (\ref{4.28}) and (\ref{3:bound})-(\ref{3:bound2}) that
\begin{equation}
\begin{split}
&\inte V(x)u_k^2dx\geq C(\delta)\int_{\big(B_{r_2}(\bar z_k)\setminus B_{r_1}(\bar z_k)\big)\cap (\Omega^\delta)^c}u_k^2dx\\
&=C(\delta)\int_{\big(B_{\frac{r_2}{\eps_k}}(0)\setminus B_{\frac{r_1}{\eps_k}}(0)\big)\cap (\Omega_{\eps_k}^\delta)^c}\bar w_k^2dx\geq C(\delta,\sigma)e^{-2(1+\sigma)\frac{r_2}{\eps_k}}\,,
\end{split}
\end{equation}
where $\Omega_{\eps_k}^\delta=\{x\in\R^2: \eps_k x+\bar z_k\in \Omega^\delta\}$.
It then follows from the Gagliardo-Nirenberg inequality (\ref{GNineq}) that
\begin{equation}\label{3:fin}
\begin{split}
e(a_k)=E_{a_k}(u_k)
&\geq
\Big(1-\frac{a_k}{a^*}\Big)\inte|\nabla u_k|^2dx+\inte V(x)u_k^2dx\\
&\geq\frac{a^*-a_k}{a^*\eps_k^2}+Ce^{-2(1+\sigma)\frac{r_2}{\eps_k}}\quad \text{as}\quad a_k\nearrow a^*\,.
\end{split}
\end{equation}
One can check that the function
\begin{equation*}
f_\gamma(s)=\frac{a^*-a}{a^*}s^2+Ce^{-\gamma s}, \quad\gamma>0,\quad s\in(0,\infty)\,,
\end{equation*}
 has a unique minimum point  $s_\gamma$   satisfying
 $$s_\gamma=\frac{1}{\gamma}\big[|\ln(a^*-a)|-\ln|\ln(a^*-a)|+C(\gamma)\big]\quad \text{as}\quad a\nearrow a^*\,.$$
Thus,
$$f_\gamma(s) \ge f_\gamma(s_\gamma)= \frac{1+o(1)}{a^*\gamma^2}(a^*-a)|\ln(a^*-a)|^2\quad \text{as}\quad a\nearrow a^*\,.$$
This estimate and (\ref{3:fin}) imply that
\begin{equation}\label{4.33}
e(a_k)\geq \frac{1+o(1)}{a^*[2(1+\sigma)r_2]^2}(a^*-a_k)\big|\ln(a^*-a_k)\big|^2\quad \text{as}\quad a_k\nearrow a^*\,.
\end{equation}
Using  (\ref{2:en}), we then have
$$(1+\sigma)r_2=(1+\sigma)(R_0+4\delta)\geq R,\quad \text{for any}\quad \sigma\in(0,1) \quad \text{and}\quad 0<\delta<\min\big\{\frac{R_{i_0}}{5},\frac{d}{5}\big\}\,. $$
This yields  that $R_0\geq R$, and thus $R_0=R_{i_0}=R$ by (\ref{3:rel}). We therefore conclude   (\ref{1:R0}).

Moreover,  by applying Lemma \ref{lem2.2},  (\ref{1:R0}) and (\ref{4.33}), we obtain that (\ref{5.1.4}) holds for the sequence $\{a_k\}$.  Since the above argument can be carried out for any sequence $\{a\}$ satisfying $a\nearrow a^*$, we also conclude that (\ref{5.1.4}) holds for all $a\nearrow a^*$.

We next show that
\begin{equation}\label{3:fin1}
\lim_{k\to\infty}\frac{\eps_k\big|\ln (a^*-a_k)\big|}{2R}=1
\end{equation}
holds (passing to a subsequence  if necessary).
In fact, if
$$0\leq\liminf_{k\to\infty}\frac{\eps_k\big|\ln (a^*-a_k)\big|}{2R}=\eta<1\,,$$
we then have
$$\frac{1}{\eps_k}\geq \frac{|\ln (a^*-a_k)|}{ 2R\eta_0}\ \, \mbox{as}\, \ k\to\infty,\quad \mbox{where}\ \,
\eta _0:=\frac{1+\eta}{2}<1.$$
We thus obtain from (\ref{3:fin}) that
$$e(a_k)\geq\frac{a^*-a_k}{a^*\eps_k^2}\geq \frac{(a^*-a_k)\big|\ln(a^*-a_k)\big|^2}{a^*[2R\eta _0]^2}\ \, \mbox{as}\, \ k\to\infty.$$
Applying Lemma \ref{lem2.2} yields that
$\eta_0\geq 1$, a contradiction. Similarly, if
$$\limsup_{k\to\infty}\frac{\eps_k\big|\ln (a^*-a_k)\big|}{2R}=\eta>1\,,$$
then $$\frac{1}{\eps_k}\leq \frac{|\ln (a^*-a_k)|}{ 2R\eta _0}\ \, \mbox{as}\ \, k\to\infty,\quad \mbox{where}\ \,
\eta _0:=\frac{1+\eta}{2}>1.$$
Since $\eta _0>1$, we are able to choose positive constants $\sigma$ and $\delta$ small enough that
$$\frac{(1+\sigma)r_2}{\eta_0 R}=\frac{(1+\sigma)(R+4\delta)}{\eta_0 R}<1\,.$$
It then follows from (\ref{3:fin}) that
$$e(a_k)\geq Ce^{-2(1+\sigma)\frac{r_2}{\eps_k}}\geq C (a^*-a_k)^\frac{(1+\sigma)r_2}{  R\eta _0}\ \mbox{as}\ k\to\infty ,$$
which however contradicts Lemma \ref{lem2.2}. We therefore conclude that  (\ref{3:fin1}) holds.

It finally follows from (\ref{4.28}) and (\ref{3:fin1}) that, for $k\to\infty$,
\begin{equation*}
\begin{split}
&\ \ \ \frac{2R}{|\ln(a^*-a_k)|} u_{a_k}\Big(\bar z_k +
\frac{2R}{|\ln(a^*-a_k)|}x\Big) \\
&=\frac{2R}{\varepsilon_k|\ln(a^*-a_k)|}\cdot\varepsilon_k u_k\Big(\varepsilon_k\cdot\frac{2R}{\varepsilon_k|\ln(a^*-a_k)|}x+\bar z_k\Big)\\
&=\frac{2R}{\varepsilon_k|\ln(a^*-a_k)|}\bar w_k\Big(\frac{2R}{\varepsilon_k|\ln(a^*-a_k)|}x\Big)
\to
 \frac { Q(x)}{\sqrt{a^*}}\quad \text{strongly in }\ H^1(\R^2),
 \end{split}
\end{equation*}
i.e., (\ref{con:a}) holds, and we therefore complete the proof of Theorem \ref{thm1.3}. \qed

%\vspace {.95cm}
\noindent {\bf Acknowledgements:} This research was  supported
    by NSFC grants 11322104, 11471331, 11271360 and 11271201.


\begin{thebibliography}{29}
%\bibitem{Adams} R. A. Adams and  J.   F. Fournier, {\it Sobolev spaces}, $2^{\rm nd}$ ed., Academic Press (2003).

%\bibitem{Anderson} M. H. Anderson, J. R. Ensher, M. R. Matthews,  C. E. Wieman and E. A. Cornell, {Observation of Bose-Einstein condensation in a dilute atomic vapor}, Science {\bf 269} (1995), 198--201.

\bibitem{Schnee} W. H. Aschbacher, J. Fr\"ohlich, G. M. Graf, K. Schnee and M. Troyer, {Symmetry breaking regime in the nonlinear Hartree equation},  J. Math. Phys. {\bf 43} (2002), 3879--3891.

\bibitem{Bao2} W. Z. Bao and Y. Y. Cai, {Mathematical theory and numerical methods for Bose-Einstein condensation}, Kinetic and Related Models {\bf 6} (2013), 1--135.

%\bibitem{Zwerger} I. Bloch, J. Dalibard and W. Zwerger, {Many-body physics with ultracold gases}, Rev. Modern Phys. {\bf 80} (2008), 885.

\bibitem{BW} T. Bartsch amd Z.-Q. Wang, {Existence and multiplicity results for some superlinear elliptic problems on $R^N$}, Comm. Partial Differential Equations  {\bf 20}  (1995),  1725--1741.


\bibitem{BO} J. Byeon and Y. Oshita, {Uniqueness of standing waves  for  nonlinear Schr\"o\-dinger equations}, Proc. Roy. Soc. Edinburgh Sect. A {\bf 138}  (2008), 975--987.

\bibitem{BWang} J. Byeon and Z. Q. Wang, {Standing waves with a critical frequency for nonlinear Schr\"o\-dinger equations}, Arch. Ration. Mech. Anal. {\bf 165}  (2002), 295--316.


\bibitem{C} T. Cazenave,  {\em Semilinear Schr\"odinger Equations}, Courant Lecture Notes in Math. {\bf 10}, Courant Institute of Mathematical Science/AMS, New York, (2003).

%\bibitem{Cooper} N. R. Cooper, { Rapidly rotating atomic gases}, Adv. Phys. {\bf 57} (2008), 539--616.

\bibitem{D} F. Dalfovo, S. Giorgini, L. P. Pitaevskii and S. Stringari, {Theory of Bose-Einstein condensation in trapped gases}, Rev. Modern Phys. {\bf 71} (1999), 463--512.

%\bibitem{Ketterle} K. B. Davis, M. O. Mewes, M. R. Andrews, N. J. van Druten, D. S. Durfee, D. M. Kurn and W. Ketterle, {\it Bose-Einstein condensation in a gas of sodium atoms}, Phys. Rev. Lett. {\bf 75} (1995), 3969--3973.

%\bibitem{DL} Y. H. Ding and S. J. Li, {Homoclinic orbits for first order Hamiltonian systems}, J. Math.
%Anal. Appl. {\bf 189} (1995), 585--601.

%\bibitem{Fetter} A. L. Fetter, {Rotating trapped Bose-Einstein condensates}, Rev. Modern Phys. {\bf 81} (2009), 647.

\bibitem{GNN} B. Gidas, W. M. Ni and L. Nirenberg, {Symmetry of positive solutions of nonlinear elliptic equations in $\R^n$}, Mathematical analysis and applications  Part A, Adv. in Math. Suppl. Stud. Vol. {\bf 7} (1981), 369--402.

\bibitem{GT}D. Gilbarg and N. S. Trudinger,  {\em Elliptic Partial Differential Equations of Second Order}, Springer, (1997).

%\bibitem{G60} E. P. Gross, {Structure of a quantized vortex in boson systems}, Nuovo Cimento {\bf 20} (1961), 454--466.

%\bibitem{G63} E. P. Gross, {Hydrodynamics of a superfluid condensate}, J. Math. Phys. {\bf 4} (1963), 195--207.


\bibitem{GS} Y. J. Guo and R. Seiringer, {On the mass concentration for Bose-Einstein condensates with attractive interactions}, Lett. Math. Phys. {\bf104} (2014), 141--156.

\bibitem{GZZ} Y. J. Guo, X. Y. Zeng and  H. S. Zhou, {Energy estimates and symmetry breaking  in attractive Bose-Einstein condensates with ring-shaped potentials}, submitted  (2014).


%\bibitem{GM} S. Gupta, K. W. Murch, K. L. Moore, T. P. Purdy and D. M. Stamper-Kurn, {Bose-Einstein condensation in a circular waveguide}, Phys. Rev. Lett. {\bf 95} (2005), 143201.

\bibitem {HL}  Q. Han and F. H. Lin,   {\em Elliptic Partial Differential Equations}, Courant Lecture Notes in Math. 1, Courant Institute of Mathematical Science/AMS, New York, 2011.


\bibitem {HM}  C. Huepe, S. Metens, G. Dewel, P. Borckmans and M. E. Brachet, {Decay rates in attractive Bose-Einstein condensates}, Phys. Rev. Lett. {\bf 82}  (1999), 1616--1619.

\bibitem {J} R. K. Jackson and M. I. Weinstein, {Geometric analysis of bifurcation and symmetry breaking in a Gross-Pitaevskii equation}, J. Stat. Phys. {\bf 116} (2004), 881--905.

\bibitem {KM} Y. Kagan, A. E. Muryshev and G. V. Shlyapnikov, {Collapse and Bose-Einstein condensation in a trapped Bose gas with nagative scattering length}, Phys. Rev. Lett. {\bf 81} (1998), 933--937.


\bibitem{KW} O. Kavian and F. B. Weissler, {Self-similar solutions of the pseudo-conformally invariant nonlinear Schr\"{o}dinger equation}, Michigan Math. J. {\bf 41} (1994), 151--173.

\bibitem {K11}  E. W. Kirr, P. G. Kevrekidis and D. E. Pelinovsky, {Symmetry-breaking bifurcation in the nonlinear Schr\"odinger equation with symmetric potentials}, Comm. Math. Phys. {\bf 308} (2011), 795--844.

\bibitem {K08} E. W. Kirr, P. G. Kevrekidis,  E. Shlizerman and  M. I.  Weinstein,  {Symmetry-breaking bifurcation in nonlinear Schr\"odinger/Gross-Pitaevskii equations}, SIAM J. Math. Anal. {\bf 40} (2008), 566--604.



\bibitem{K} M. K. Kwong, {  Uniqueness of positive solutions of $\Delta u-u+u^p=0$ in $\R^N$}, Arch. Rational Mech. Anal. {\bf 105} (1989), 243--266.

%\bibitem{L} H. A. Levine, {\em An estimate for the best constant in a Sobolev inequality involving three integral norms}, Ann. Math. Pura Appl. {\bf 124}, 181--197 (1980).

\bibitem{Li} Y. Li and W. M. Ni, { Radial symmetry of positive solutions of nonlinear elliptic equations in $\R^n$}, Comm. Partial Differential Equations {\bf 18}  (1993), 1043--1054.

%\bibitem{Lieb} E. H. Lieb and  M. Loss, {\em Analysis}, Graduate Studies in Math. {\bf  14}, Amer. Math. Soc., Providence, RI, second edition (2001).

%\bibitem{LS} E. H. Lieb and R. Seiringer, {Proof of Bose-Einstein  condensation for dilute trapped gases}, Phys. Rev. Lett. {\bf 88} (2002), 170409-1--4.

%\bibitem {LSS} E. H. Lieb, R. Seiringer, J. P. Solovej and J. Yngvason, {\em The Mathematics of the Bose Gas and its Condensation}, Oberwolfach Seminars {\bf 34}, Birkh\"auser Verlag, Basel (2005).

%\bibitem {LSY} E. H. Lieb, R. Seiringer and J. Yngvason, {Bosons in a trap: A rigorous derivation of the Gross-Pitaevskii energy functional}, Phys. Rev. A {\bf 61} (2000), 043602-1--13.

%\bibitem{LSY2d} E. H. Lieb, R. Seiringer and J. Yngvason,  {A rigorous derivation of the Gross-Pitaevskii energy functional for a two-dimensional Bose gas}, Comm. Math. Phys. {\bf 224} (2001), 17--31.

%\bibitem{l1} P. L. Lions, {The concentration-compactness principle in the caclulus of variations. The locally compact case I}, Ann. Inst. H. Poincar\'{e} Anal. Non Lin\'{e}aire {\bf 1} (1984), 109--145.

%\bibitem{l2} P. L. Lions, {The concentration-compactness principle in the caclulus of variations. The locally compact case II}, Ann. Inst. H. Poincar\'{e} Anal. Non Lin\'{e}aire {\bf 1}  (1984), 223--283.

\bibitem{LW} G. Z. Lu and J. C. Wei, {On nonlinear schr\"odinger equations with totally degenerate potentials}, C. R. Acad. Sci. Paris. {\bf 326}  (1998), 691--696.

\bibitem{mcleod} K. McLeod and J. Serrin, {Uniqueness of positive radial solutions of $\Delta u +f (u)=0$ in $\R^n$}, Arch. Rational Mech. Anal. {\bf 99}  (1987), 115--145.

\bibitem {NT} W.-M. Ni and I. Takagi, {On the shape of least-energy solutions to a semilinear Neumann problem}, Comm. Pure Appl. Math. {\bf 44} (1991), 819--851.

\bibitem {NW} W.-M. Ni and J. C. Wei, {On the location and profile of spike-layer solutions to singularly perturbed  semilinear Dirichlet problems}, Comm. Pure Appl. Math. {\bf 48} (1995), 731--768.



\bibitem {P} L. P. Pitaevskii, {Vortex lines in an imperfect Bose gas}, Sov. Phys. JETP. {\bf 13}  (1961), 451--454.

\bibitem {RS} M. Reed and B. Simon, {\em Methods of Modern Mathematical Physics. IV. Analysis of Operators}, Academic Press, New York-London, 1978.

%\bibitem {S} J. Smyrnakis, S. Bargi, G. M. Kavoulakis, M. Magiropoulos, K. Karkkainen and S. M. Reimann, {Mixtures of Bose gases confined in a ring potential}, Phys. Rev. Lett. {\bf 103} (2009), 100404.


\bibitem {Wang} X. F. Wang, {On concentration of positive bound states of nonlinear Schr\"odinger equations}, Comm. Math. Phys. {\bf 153} (1993),  229--244.

\bibitem {WangZ} Z.-Q. Wang, {Existence and symmetry of multi-bump solutions for nonlinear Schrödinger equations}, J. Differential Equations  {\bf 159}  (1999), 102--137.

\bibitem {W} M. I. Weinstein, {Nonlinear Schr\"odinger equations and sharp interpolations estimates}, Comm. Math. Phys. {\bf 87} (1983), 567--576.

%\bibitem {Z} J. Zhang, {Stability of attractive Bose-Einstein condensates}, J. Stat. Phys. {\bf 101}  (2000), 731--746.









\end{thebibliography}
\end{document}